\documentclass[12pt,english,twoside]{article}
\usepackage[T1]{fontenc}
\usepackage[latin1]{inputenc}
\usepackage{babel}

\usepackage{geometry}
\usepackage{amsmath}

\usepackage{times}

\usepackage{amssymb}
\usepackage{verbatim}

\makeatletter

\newtheorem{thr}{Theorem}[section]
\newtheorem{lm}{Lemma}[section]

\newtheorem{prop}{Proposition}[section]
\newtheorem{dfn}{Definition}[section]

\newcommand{\ri}{\rightarrow}
\newcommand{\eps}{\varepsilon}
\newcommand{\sn}{{\scriptscriptstyle N}}
\newcommand{\uds}{U_N^+}
\newcommand{\E}{\mathsf{E}}
\renewcommand{\P}{\mathsf{P}}
\newcommand{\p}{p}
\newcommand{\py}{\pi^Y}
\newcommand{\si}{\sigma}

\newcommand{\La}{\si }
\newcommand{\lmax}{\La _{N}^{*}}
\newcommand{\lmin}{\La _{*,N}}
\newcommand{\mxi}{\xi}

\newcommand{\strtop}{\mbox{\itshape{s.t.}}}
\newcommand{\ds}{\displaystyle}

\newcommand{\R}{{\textbf{R}}}
\newcommand{\uc}{u^{\circ }}
\newcommand{\g}{\gamma}

\newcommand{\ke}{\varkappa _{\varepsilon }}
\newcommand{\supp}{\mbox{supp}\,}

\renewcommand{\@evenfoot}{\null{\footnotesize 26/01/2005}\hfill\null}
\renewcommand{\@oddfoot}{\null}
\renewcommand{\@oddhead}{{\small Asymptotic analysis of a stochastic model
for  parallel computations}\null\hfill\thepage}
\renewcommand{\@evenhead}{\thepage\null\hfill{\small A.~Manita, V.~Shcherbakov}}

\makeatother
\begin{document}

\title{Asymptotic analysis of a particle system with mean-field interaction}

\author{Anatoli Manita \thanks{Supported by
  Russian Foundation of Basic Research (RFBR grant
 02-01-00945).},\\
{\small Faculty of Mathematics and
Mechanics,}\\
{\small Moscow State University, 119992, Moscow, Russia.}\\
{\small E-mail:~manita@mech.math.msu.su}\and
Vadim Shcherbakov\thanks{Supported by Russian Foundation of Basic Research
(RFBR grant  01-01-00275). On leave from Laboratory of Large Random Systems, Faculty of
Mathematics and Mechanics, Moscow State University, 119992, Moscow, Russia.},\\
{\small CWI, Postbus 94079, 1090 GB,}\\
{\small Amsterdam, The Netherlands}\\
{\small E-mail:~V.Shcherbakov@cwi.nl}}
\date{}

\maketitle

\vspace*{-5ex}
\begin{abstract}

We study  a system of~$N$ interacting particles on~$\mathbf{Z}$.
The~stochastic dynamics consists of two components: a free motion of
each particle (independent random walks) and a~pair-wise interaction
between particles. The~interaction belongs to the class of {\it
mean-field\/} interactions and models a {\it rollback
synchronization\/} in asynchronous networks of processors for a
distributed simulation. First of all we study an  empirical measure
generated by the particle configuration on $\mathbf{R}$. We~prove that
if space, time and a parameter of the interaction are appropriately
scaled (hydrodynamical scale), then the empirical measure converges
weakly  to a deterministic limit as $N$ goes to infinity. The~limit
process is defined as a weak solution of some partial differential
equation. We~also study the~long time evolution
of the particle system with fixed number of particles. The~Markov
chain formed by individual positions of the particles is not
ergodic. Never\-theless it is possible to introduce {\it relative}
coordinates and prove that the~new Markov chain is ergodic while
the system as a whole moves with an asymptotically constant mean
speed which differs from the mean drift of the free particle motion.

\end{abstract}

\section{Introduction}

\label{intro}

We study an interacting particle system which  models a set of
processors performing parallel simulations. The system can be
described as follows. Consider $N\geq 2$ particles moving
in~${\textbf {Z}}$.  Let $x_{i}(t)$ be the position at time $t$ of
the $i-$th particle, $1\leq i\leq N$. Each particle has three
clocks. The first, the second and the third clock, attached to the
$i-$th particle, ring at the moments of time given by a mutually
independent Poisson processes $\Pi _{i,\alpha },\, \Pi _{i,\beta }$
and $\Pi _{i,\mu _{N}}$ with intensities $\alpha $, $\beta $ and
$\mu _{N}$ correspondingly. These triples of Poisson processes for
different indexes are also independent. Consider a particle with
index $i$. If the first attached clock rings, then the particle
jumps to the nearest right site: $x_{i}\ri x_{i}+1$, if the second
attached clock rings, then the particle jumps to the nearest left
site: $x_{i}\ri x_{i}- 1$. At moments when the third attached clock
rings a particle with index $j$ is chosen with probability $1/N$ and
if $x_{i}>x_{j}$, then the $i-$th particle is relocated: $x_{i}\ri
x_{j}$. It is supposed that all  these changes occur immediately.

The  type of the interaction between the particles is motivated by
studying of probabilistic models in the theory of parallel
simulations in computer science (\cite{Mitra}, \cite{T1,T3} and
\cite{GrMP,MM-Time-Sync}).  The main peculiarity of the models is
that a group of processors performing a large-scale simulation is
considered and each processor does a specific part of the task. The
processors share data doing simulations therefore their activity
must be synchronized. In practice, this synchronization is achieved
by applying a so-called {\it rollback\/} procedure which is based
on a~massive message exchange between different processors
(see~\cite[Sect.~1.4 and Ch.~8]{BertTsit}). One says that $x_i(t)$
is a {\it local time\/} of the $i$-th processor while $t$ is a real
(absolute) time. If~we interpret the variable $x_i(t)$ as an amount
of job done by the processor~$i$ till the time moment~$t$, then the
interaction described above imitates this  synchronization
procedure. Note that from a point of view of general stochastic
particle systems the interaction between  the particles is
essentially non-local.

We are interested in the analysis of asymptotic behaviour of this
particle system. First of all we consider  the situation as  the
number of particles goes to infinity. For every finite $N$ and $t$
we can define an  empirical  measure generated by the particle
configuration. It~is a point measure with atoms at integer points.
An atom at a point $k$ equals to a proportion of particles with
coordinate $k$ at time $t$. It is convenient in our case to consider
an empirical tail function corresponding to the measure. It means
that we consider 
$
\xi _{x,N}(t)={\ds \frac{1}{N}}\, 
\sum \limits _{i=1}^{N}\mathbf{1}\left(x_{i}(t)\geq x\right)
$
 the proportion of particles having
coordinates not less than $x\in \mathbf{R}$. The problem is to find
an appropriate time scale $t_{N}$ and a sequence of interaction
parameters $\mu_{N}$ to obtain a non-trivial limit dynamics of the
process $\xi_{N, [xN]}(t_{N})$ as $N\ri \infty $. The cases
$\alpha\neq \beta$ and $\alpha=\beta$ require different  scaling of
time and the interaction constant $\mu_{N}$. We prove that there
exist non-trivial limit deterministic processes in both cases as $N$
goes to infinity if we rescale time and the interaction constant as
$t_{N}=tN,\,\mu_{N}=\mu/N$ in the first case and as
$t_{N}=tN^{2},\,\mu_{N}=\mu/N^{2}$ in the second case respectively.
The processes are defined as weak solutions of some partial
differential equations (PDE). 
 It should be noted that
the PDE relating to the zero drift
situation is a famous \emph{Kolmogorov--Petrowski--Piscounov}-equation
(KPP-equation, \cite{KoPP}). This result was announced in
\cite{MaSch}.

Another issue we address in the paper is studying of the long time
evolution of the particle system with fixed number of particles. It
is easy to see that the Markov chain $x(t)=\{x_{i}(t),\,i=1,\ldots
,N\},\,t\geq 0,$ is not ergodic. Nevertheless the particle system
possesses some relative stability. We introduce new coordinates
$y_{i}(t)=x_{i}(t)-\min_{j}x_{j}(t),\,i=1,\ldots ,N$, and
prove that the countable Markov chain $y(t)=\{y_{i}(t),\,i=1,\ldots
,N\},\,t\geq 0,$ is ergodic and converges exponentially fast to its
stationary distribution. Therefore the system of stochastic
interacting particles possesses some relative stability.  We show
also that the center of mass of the system moves with an
asymptotically constant speed. It appears  that due to the
interaction between the particles this speed differs from the mean
drift of the free particle motion.

It should be noted that the choice  of the interaction  may vary
depending on a situation. Various modifications
of the model can be considered and similar results can be obtained
using the same methods. We have chosen the described model
just for the sake of concreteness.

Probabilistic models for parallel computation  considered before by 
other authors. The paper~\cite{Mitra} deals with a model consisting of 
two interacting processors ($N=2$). It contains a rigorous study of the 
long-time behavior of the system and formulae for some performance 
characteristics. Unfortunately, there are not too many mathematical 
results about multi-processor models 
(\cite{MadWalMes,GuAkFu,AkChFuSer,T1,T3}). Usually mathematical 
components of these papers have a form of preparatory considerations 
before some large numerical simulation. The paper~\cite{GrMP} is of 
special interest because it rigorously studies a behavior of some model 
of parallel computation with $N$ processor units in the limit 
$N\rightarrow\infty$. A stochastic dynamics of~\cite{GrMP} is different 
from the dynamics studied in the present paper and main results 
of~\cite{GrMP} concern a so-called thermodynamical limit. The authors 
prove that in the limit the evolution of the system can be described by 
some integro-differential equation. In the present study we propose a 
model which dynamics is easy from the point of view of numerical 
simulations and, at the same time,  provides  us with a new 
probabilistic interpretation of some important PDEs including 
the classical KPP-equation.

The paper is organised as follows. We formally define the particle 
system, introduce some notation and formulate the main results in 
Section \ref{model}. Sections \ref{proof} and \ref{stability} contain 
the proofs of the main results. In Section \ref{travel} we discuss 
solutions of the limiting equations.

\paragraph{Acknowledgments.}

We are thankful to Dr.\ T.~Voznesenskaya (Faculty of Computational 
Mathematics and Cybernetics, Moscow University) who first introduced us 
to stochastic algorithms for parallel computations. The authors would 
like to thank Prof.\ V.~Bogachev (Faculty of Mechanics and Mathematics, 
Moscow University) for the  helpful  discussions on convergence of 
measures on general topological vector spaces and for the suggested 
references. We are  also grateful to Prof.\ V.~Malyshev for his warm 
encouragement and valuable comments on the present  manuscript.

\section{The model and main results}
\label{model}

Formally, the process $x(t)=\{x_{i}(t),i=1,\ldots ,N\}$, describing
positions of the particles, is a continuous time countable Markov
chain taking values in ${\textbf {Z}}^{N}$ and having  the following
generator
\begin{multline}
G_{N}g(x)=\sum\limits _{i=1}^{N} \left(
\alpha \left(g\left(x+e_{i}^{\left(
\sn \right)}\right)-g(x)\right)+
\beta \left(g\left(x-e_{i}^{\left(\sn \right)}\right)-g(x)\right)
\right)+\\
{}+\sum \limits _{i=1}^{N}\sum \limits _{j\neq i}\left(g\left(x-e_{i}^{\left(
\sn \right)}
\left(x_{i}-x_{j}\right)\right)-g(x)\right)I_{\{x_{i}>x_{j}\}}
\frac{\mu_{N}}{N}\, ,
\end{multline}
 where $x=(x_{1},\ldots ,x_{N})\in {\textbf {Z}}^{N}$,
 $g:\;{\textbf {Z}}^{N}\rightarrow \textbf {R}$ is a bounded function,
 $e_{i}^{\left(\sn \right)}$
is a $N$-dimensional vector with all zero components except $i-$th
which equals to $1$, $I_{\{x_{i}>x_{j}\}}$ is an indicator of the
set $\{x_{i}>x_{j}\}$.

Define
\begin{equation}
\xi _{N,k}(t)=\frac{1}{N}\, \sum _{i=1}^{N}I_{\{x_{i}(t)\geq k\}},\quad k\in
\mathbf{Z}.
\label{ksi}
\end{equation}
The process $\xi _{N}(t)=\{\xi _{N,k}(t),\, k\in \mathbf{Z}\}$ is
a Markov one with a state space $H_{N}$ the set of all non-negative
and nonincreasing sequences $z=\{z_{k},k\in {\textbf {Z}}\}$ such
that $z_{k}\in \{l/N,\,l=0,1,\ldots,N\}$ for any $k\in {\textbf {Z}}$
and
\[
\lim _{k\rightarrow -\infty }z_{k}=1,\qquad \lim _{k\rightarrow +\infty
 }z_{k}=0.
\]

The generator of the process $\xi _{N}(t)$ is given by the following formula
\begin{align}
L_{N}f(z) & =N\sum \limits _{k}((f(z+e_{k}/N)-f(z))\alpha
 (z_{k-1}-z_{k})+(f(z-e_{k}/N)-f(z))\beta (z_{k}-z_{k+1}))\label{L2}\\
 & +N\mu _{N}\sum \limits _{l<k}(f(z-(e_{l+1}+\ldots
 +e_{k})/N)-f(z))(z_{k}-z_{k+1})(z_{l}-z_{l+1})\, ,\nonumber
\end{align}
where $e_{i},i\in {\textbf {Z}}$ is an infinite dimensional vector
with all zero components except $i-$th which equals to $1$,
$f:\;H_{N}\rightarrow \textbf {R}$ is a bounded function.

\subsection{Hydrodynamical behavior of the particle system}
\label{hydro}

Denote $\zeta_{N,x}(t)=\xi_{N,[Nx]}(t),\, \, x\in \mathbf{R}$. The
process $\zeta _{N}(t)$ takes values in $H=H(\mathbf{R})$ the set of
all non-negative right continuous with left limits nonincreasing
functions having the following limits
\[
\lim _{x\ri -\infty }\psi (x)=1,\, \lim _{x\ri \infty }\psi (x)=0.
\]
Denote by $S(\textbf{R})$ the Schwartz space of infinitely differentiable functions
such that for all $m,n\in\textbf{Z}_+$
$$
\|f\|_{m,n}=\sup_{x\in\textbf{R}}|x^m f^{(n)}(x)|<\infty .
$$
Recall that $S(\textbf{R})$ equipped with a natural topology given by seminorms
$\|\cdot\|_{m,n}$ is a Frechet space~(\cite{RS}).

Define for every $h\in H$ a functional
\[
(h,f)=\int \limits _{\textbf{R}}h(x)f(x)dx,\,\, f\in S(\textbf{R}),
\]
on  the Schwartz space $S(\textbf{R})$.
The following bound yields that for each $h\in H$ $(h,\cdot)$ is a~continuous
linear functional on $S(\textbf{R})$
$$
|(h,f)|\leq \int \limits _{\textbf{R}}|f(x)|\frac{1+x^2}{1+x^2}dx\leq
\pi(\|f\|_{\infty}+\|x^{2}f\|_{\infty})\equiv
\pi(\|f\|_{0,0}+\|x^{2}f\|_{2,0}),
$$
where $\|\cdot\|_{\infty}$ is the supremum norm. Thus the set of
functions $H(\textbf{R})$ is naturally embedded into the space of
all continuous linear functionals on $S({\textbf{R}})$, namely into
the space $S'(\textbf{R})$ of tempered distributions. We will
interpret $\zeta _{N}(t)$ as a stochastic process taking its values
in the space $S'(\textbf{R})$.

There are two reasons for embedding $H(\textbf{R})$ into 
$S'(\textbf{R})$ and considering the $S'(\textbf{R})$-valued processes. 
The first reason is that due to some nice topological properties of 
$S'(\textbf{R})$ we can use in Section~\ref{proof} many powerful results 
from the theory of weak convergence of probability distributions on topological 
vector fields. And, secondly, the choice of $S'(\textbf{R})$ as a state 
space is convenient from the point of view of possible future study of 
stochastic fluctuation fields around the deterministic limits obtained 
in our main theorem~\ref{hydro}.

In the sequel we mainly deal with the strong topology ($\strtop$) on 
$S'(\textbf{R})$ (see Section~\ref{top}). From now we fix some $T>0$ and 
consider $\zeta _{N}$ as a random element in a Skorokhod space $D([0,T], 
S'(\textbf{R}))$ of all mappings of $[0,T]$ to $(S'(\textbf{R}),\strtop)$ 
that are right continuous and have left-hand limits in the strong 
topology on $S'(\textbf{R})$. Note that $(S'(\textbf{R}),\strtop)$ is not 
a metrisable topological space therefore it is not evident how to 
define the  Skorokhod topology on the space $D([0,T], S'(\textbf{R}))$. 
To do this we follow~\cite{Mitoma} and refer to Section~\ref{top}.

Now we are able to consider  probability distributions  of the
processes $\left(\zeta _{N}(tN^{a}), t\in [0,T]\right)$, $a=1,2,$ as
probability measures on a measurable space $\left(D([0,T],
S'(\textbf{R})),\mathcal{B}_{D([0,T], S'(\textbf{R}))}\right)$ where
$\mathcal{B}_{D([0,T], S'(\textbf{R}))}$ is a corresponding Borel
$\sigma$-algebra.
It was proved in~\cite{Jakub} that 
$\mathcal{B}_{D([0,T], S'(\textbf{R}))}=\mathcal{C}_{D([0,T], S'(\textbf{R}))}$,
where $\mathcal{C}_{D([0,T], S'(\textbf{R}))}$ is a $\sigma$-algebra of cylindrical subsets.

Consider two following Cauchy problems
\begin{align}
u_{t}(t,x) & =-\lambda u_{x}(t,x)+\mu (u^{2}(t,x)-u(t,x))\,
 ,\label{eq:ur-asym}\\
 u(0,x) & =\psi (x)\nonumber
\end{align}
 and
\begin{align}
u_{t}(t,x) & =\gamma u_{xx}(t,x)+\mu (u^{2}(t,x)-u(t,x))\, ,\label{eq:ur-KPP}\\
u(0,x) & =\psi (x)\nonumber
\end{align}
where $u_{t},u_{x}$ and $u_{xx}$ are first and second derivatives of $u$ with
respect to $t$ and $x$. Notice that the equation~(\ref{eq:ur-KPP}) is a
particular case of the famous \emph{Kolmogorov--Petrowski--Piscounov}-equation
(KPP-equation, \cite{KoPP}). We will deal with   weak solutions of the equations
(\ref{eq:ur-asym}) and (\ref{eq:ur-KPP}) in the sense of Definition \ref{weak}.

Fix $T>0$ and denote by $C_{0,T}^{\infty}=C_{0}^{\infty}([0,T]\times
{\textbf {R}})$ the space of infinitely differentiable functions
with finite support and equal to zero for $t=T$.

\begin{dfn}
\label{weak}
\begin{description}
\item[(i)]
The bounded measurable function \(u(t,x)\) is called a weak (or generalized)
solution of the Cauchy problem  (\ref{eq:ur-asym})
in the region $[0,T]\times \textbf{R}$,
if the following integral equation holds for any function
\(f \in C_{0,T}^{\infty}\)
\begin{multline*}
\int\limits_{0}^{T}\int\limits_{\textbf{R}}u(t,x)(f_{t}(t,x)
+\lambda f_{x}(t,x))
+\mu u(t,x)(1-u(t,x))f(t,x))dxdt\\
+\int\limits_{\textbf{R}}u(0,x)f(0,x)dx=0
\end{multline*} 
\item[(ii)]
The bounded measurable function \(u(t,x)\) is called a weak (or generalized)
solution of the Cauchy problem  (\ref{eq:ur-KPP})
in the region $[0,T]\times \textbf{R}$,
if the following integral equation holds for any function
\(f \in C_{0,T}^{\infty}\)
\begin{multline*}
\int\limits_{0}^{T}\int\limits_{\textbf{R}}u(t,x)(f_{t}(t,x)
+\gamma f_{xx}(t,x))
+\mu u(t,x)(1-u(t,x))f(t,x))dxdt\\
+\int\limits_{\textbf{R}}u(0,x)f(0,x)dx=0
\end{multline*}
\end{description}
\end{dfn}

In Subsection~\ref{edinstv} we will show that the both of Cauchy 
problems (\ref{eq:ur-asym}) and (\ref{eq:ur-KPP}) have \emph{unique weak 
solutions} in the sense of Definition~\ref{weak}. Here we want just to 
mention that this problem is not trivial. Indeed, the equation 
(\ref{eq:ur-asym}) is an example of a quasilinear first order partial 
differential equation. It is known that in a general case such type of 
equations might have more than one weak solution and  it is only 
possible to guarantee uniqueness of the solution which satisfies  to the 
so-called entropy condition. The most general form of this condition was 
introduced by Kruzhkov in \cite{Kruzhkov}, where he also proved his 
famous uniqueness theorem. Fortunately, in our particular case of the 
equation (\ref{eq:ur-asym}) the situation is quite simple due to 
simplicity of characteristics,  they  are given by the straight lines 
$x(t)=\lambda t+ C$, do not intersect with each other and do not produce 
the shocks. Detailed discussions of the problem of uniqueness for 
equations~(\ref{eq:ur-asym}) and (\ref{eq:ur-KPP}) are presented in 
Subsection~\ref{edinstv}.

The first theorem  we are formulating describes
the evolution of the system at the hydrodynamical scale.
\begin{thr}
\label{thydro}
Assume that an initial particle  configuration
\(\xi_{N}(0)=\{\xi_{N,k}(0),\,k\in {\textbf{Z}}\}\)
is such that for any function \(f\in S(\textbf{R})\)
\begin{equation}
\label{initial}
\lim\limits_{N\ri\infty}\frac{1}{N}\sum\limits_{k}\xi_{N,k}(0)f(k/N)=
\int\limits_{{\textbf{R}}}\psi(x)f(x)dx,
\end{equation}
where \(\psi \in H(\textbf{R})\).
\begin{description}
 \item[(i)]
If \(\alpha-\beta=\lambda\neq 0\) and \(\mu_{N}=\mu/N\), then the
sequence $\{Q_{N,\lambda}^{\left(T\right)}\}_{N=2}^{\infty}$ of
probability distributions of random processes $\{\zeta_{N}(tN),\,
t\in [0,T]\}_{N=2}^{\infty}$ converges weakly  as \(N\ri \infty\) to
the probability measure $Q_{\lambda}^{\left(T\right)}$ on $D([0,T],
S'(\textbf{R}))$ supported by a trajectory $u(t,x)$,  which is a
unique weak solution of the equation (\ref{eq:ur-asym}) with the
initial condition \(u(0,x)=\psi(x)\) and as a function of $x$
$u(t,\cdot)\in H(\textbf{R}),$ for any $t\geq 0$.

\item[(ii)]
If $\alpha=\beta=\gamma>0\), \(\mu_{N}=\mu/N^{2}$, then the
sequence $\{Q_{N,\gamma}^{\left(T\right)}\}_{N=2}^{\infty}$ of
probability distributions  of random processes
$\{\zeta_{N}(tN^{2}),\, t\in [0,T]\}_{N=2}^{\infty}$ converges
weakly as \(N\ri \infty\) to the probability measure
$Q_{\gamma}^{\left(T\right)}$ on $D([0,T], S'(\textbf{R}))$
supported by a  trajectory $u(t,x)$,  which is a  unique weak
solution of the equation (\ref{eq:ur-KPP}) with the initial
condition \(u(0,x)=\psi(x)\) and as a function of $x$
$u(t,\cdot)\in H(\textbf{R}),$ for any $t\geq 0$.

\end{description}
\end{thr}

\subsection{Long time behavior of the particle system with the fixed
number of particles}
\label{connection}

The number of particles is fixed in this section.
Consider the following stochastic process $y(t)=(y_{1}(t),\ldots ,y_{N}(t))$,
where
$$y_{i}(t)=x_{i}(t)-\min\limits_{j}x_{j}(t).$$
Note that $x_k-x_l=y_k-y_l$ for any pair $k,l$.

It is easy to see that $y(t)$ is a continuous time Markov chain on
the state space
$$\Gamma =\bigcup _{k}\Gamma _{k}\subset \mathbf{Z}_{+}^{N},$$
where $\Gamma _{k}:=\left\{ \left(z_{1},\ldots ,z_{k-1},0,z_{k+1},
\ldots,z_{N}\right):\, z_{j}\in \mathbf{Z}_{+}\right\}$.
\begin{thr}
\label{ergodic}
The  Markov chain
 $\left(y(t),t\geq 0\right)$
is ergodic and converges exponentially fast to its stationary distribution
$$
\sum_{y\in\Gamma} |P(y(t)=y)-\pi(y)| \leq C_{1}\exp(-C_{2}t)
$$
uniformly in initial distributions of $y(0)$.
\end{thr}

\section{Proof of Theorem \ref{thydro}}
\label{proof}

\subsection{Plan of the proof}

The proof of the convergence uses the next well-known general idea (see, 
for example, \cite[\S~5]{SmFom}). Let $\{a_n\}$ be a sequence in some 
Hausdorff topological space and assume that $\{a_n\}$  satisfies to the 
following two properties:
(a) for any subsequence of $\{a_n\}$ there is a converging subsequence
(this property is called a \emph{sequential compactness});
(b) $\{a_n\}$  contains at most one limit point.
Then the sequence  $\{a_n\}$ has a limit.

In our situation the role of  $\{a_n\}$ is played by the sequences
$\{Q_{N,\lambda}^{\left(T\right)}\}_{N=2}^{\infty}$
and  $\{Q_{N, \gamma}^{\left(T\right)}\}_{N=2}^{\infty}$.
Our proof consists of the following steps.

\emph{Step 1}.
We fix an arbitrary  $T>0$ and prove that  the sequences of
probability measures $\{Q_{N,\lambda}^{\left(T\right)}\}_{N=2}^{\infty}$
and  $\{Q_{N, \gamma}^{\left(T\right)}\}_{N=2}^{\infty}$
 are  tight. We  use the Mitoma theorem (\cite{Mitoma}) and apply martingale techniques
widely used in  the theory of hydrodynamical limits of interacting particle systems
(\cite{PreDeMas,KipLan}).

It is important to note  that if a topological space $\mathcal{V}$ 
is not metrisable then, 
generally speaking, the tightness of a family of distributions on 
$\mathcal{V}$ does not imply a sequential compactness (see, for example, 
\cite[V.~2, \S~8.6]{VBogach}). So, in general, the above property~(a) 
does not follow directly from the step~1. But in
our concrete case $\mathcal{V}=D([0,T], S'(\textbf{R}))$ 
we can proceed as follows. It was shown in~\cite{Jakub} that any compact 
subset of $D([0,T], S'(\textbf{R}))$ is metrisable. Due to this property 
we can apply the theorem from~\cite[Th.~2, \S~5]{SmFom} which states that 
(under assumption of metrisability of compact subsets) the tightness of 
a family of measures implies its sequential compactness. All this 
justifies the next step.

\emph{Step 2}.
We show that a measure that is a limit of some subsequence of the 
sequence $\{Q_{N,\lambda}^{\left(T\right)}\}_{N=2}^{\infty}$ (or 
$\{Q_{N, \gamma}^{\left(T\right)}\}_{N=2}^{\infty}$) is supported by the 
weak solutions of the partial differential equation  (\ref{eq:ur-asym}) 
(or, correspondingly, (\ref{eq:ur-KPP})\,). Then we note that each of 
the equations  (\ref{eq:ur-asym}) and (\ref{eq:ur-KPP}) has a unique 
weak solution (Subsection~\ref{edinstv}). This gives the above property (b).

\subsection{Technical lemmas}

We start with  some bounds which will be used throughout  the proof.
Denote
\[
R_{f}(z)=\frac{1}{N}\sum _{k}f(k/N)z_{k},
\]
for $z\in H_{N}$ and $f\in S(\R)$.

\begin{lm}
\label{lb1}
\begin{description}
 \item[(i)] If \(\alpha\neq \beta\) and \(\mu_{N}=\mu/N\), then
 for  any  $z\in H_{N}$
\begin{equation}
\label{b1}
\left|L_{N}R_{f}(z)\right|\leq \frac{C}{N},
\end{equation}
and
\begin{equation}
\label{b2}
N\left(L_{N}R_{f}^{2}(z)-2R_{f}(z)L_{N}R_{f}(z)\right)=
O\left(\frac{1}{N}\right).
\end{equation}
\item[(ii)] If \(\alpha=\beta\) and \(\mu_{N}=\mu/N^2\), then
for  any  $z\in H_{N}$
\begin{equation}
\label{b12}
\left|L_{N}R_{f}(z)\right|\leq \frac{C}{N^2},
\end{equation}
and
\begin{equation}
\label{b22}
N^2\left(L_{N}R_{f}^{2}(z)-2R_{f}(z)L_{N}R_{f}(z)\right)=
O\left(\frac{1}{N^2}\right).
\end{equation}
\end{description}
In both cases $C=C(f,\alpha ,\beta ,\mu )$.
\end{lm}
\paragraph{Proof of Lemma \ref{lb1}.}
We will prove the bounds (\ref{b1}) and (\ref{b2}), the other ones can be
proved similarly.
We start with the bound (\ref{b1}).
Using the equations
\begin{align*}
R_{f}(z+e_{k}/N)-R_{f}(z)&  =\frac{f(k/N)}{N^{2}},\\
R_{f}(z-e_{k}/N)-R_{f}(z)& =-\frac{f(k/N)}{N^{2}}
\end{align*}
we get that for every $z\in H_{N}$
\begin{align*}
L_{N}R_{f}(z)= & \frac{1}{N}\sum \limits _{k}z_{k}(\beta f((k-1)/N)-(\alpha
 +\beta )f(k/N)+\alpha f((k+1)/N)\\
 & -\frac{\mu }{N^{2}}\sum\limits_{k}f(k/N)z_{k}(1-z_{k}).
\end{align*}

For any function $f\in S\left(\mathbf{R}\right)$ consider
its upper Darboux sum
\begin{equation}    \label{e-uds}
\uds (f)=\, \frac{1}{N}\, \sum _{k\in \mathbf{Z}}\, \max_{y\in
 \left[k/N,(k+1)/N\right]}f\left(y\right).
\end{equation}
 Since $\uds (f)\rightarrow \int f(x)\, dx$ as $N\rightarrow \infty $,
the sequence $\left\{ \uds (f)\right\} _{N=1}^{\infty }$ is bounded
in $N$ for any fixed $f$. We have uniformly in~$z\in H_{N}$
\begin{eqnarray*}
\left|L_{N}R_{f}(z)\right| & \leq & \frac{1}{N}\left(
\sum_{k}\left|\alpha (f\left((k+1)/N\right)-f\left(k/N\right))+
 \beta (f\left(k/N\right)-f\left((k-1)/N\right) )\right|\right)\\
 & + & \frac{\mu }{N}\left(\frac{1}{N}\sum _{k}\left|f(k/N)\right|\right)\\
 & \leq  & \frac{1}{N}\left(|\alpha -\beta |\uds
 \left(\left|f_{x}\right|\right)+\mu \uds \left(\left|f\right|\right)\right),
\end{eqnarray*}
where $f_{x}=df(x)/dx$. So the bound (\ref{b1}) is proved.

Let us prove the bound (\ref{b2}).
Note that $L_{N}=L_{N}^{\left(0\right)}+L_{N}^{\left(1\right)}$,
where
\begin{multline*}
L_{N}^{\left(0\right)}f(z)=N\sum _{k}(f(z+e_{k}/N)-f(z))\alpha
 (z_{k-1}-z_{k})+\\
N\sum _{k}(f(z-e_{k}/N)-f(z))\beta (z_{k}-z_{k+1})
\end{multline*}
and
\[
L_{N}^{\left(1\right)}f(z)=\mu\sum _{l<k}(f(z-(e_{l+1}+\ldots
 +e_{k})/N)-f(z))(z_{k}-z_{k+1})(z_{l}-z_{l+1}).
\]
Using the equations
\begin{align*}
R_{f}^{2}(z+e_{k}/N)-R_{f}^{2}(z) &
 =\left(2R_{f}(z)+\frac{f(k/N)}{N^{2}}\right)\frac{f(k/N)}{N^{2}},\\
R_{f}^{2}(z-e_{k}/N)-R_{f}^{2}(z) &
 =-\left(2R_{f}(z)-\frac{f(k/N)}{N^{2}}\right)\frac{f(k/N)}{N^{2}},
\end{align*}
one can obtain that  for any $z\in H_{N}$
\begin{align}
L_{N}^{\left(0\right)}R_{f}^{2}(z) & =\alpha N\sum \limits
 _{k}(R_{f}^{2}(z+e_{k}/N)-R_{f}^{2}(z))(z_{k-1}-z_{k})\nonumber \\
 & +\beta N\sum \limits
 _{k}(R_{f}^{2}(z-e_{k}/N)-R_{f}^{2}(z))(z_{k}-z_{k+1})\nonumber \\
 & =2R_{f}(z)\frac{\alpha }{N}\sum \limits
 _{k}f(k/N)(z_{k-1}-z_{k})+\frac{\alpha }{N^{3}}\sum \limits
 _{k}f^{2}(k/N)(z_{k-1}-z_{k})\nonumber \\
 & -2R_{f}(z)\frac{\beta }{N}\sum \limits _{k}f(k/N)(z_{k}-z_{k+1})+\frac{\beta
 }{N^{3}}\sum \limits _{k}f^{2}(k/N)(z_{k}-z_{k+1})\nonumber \\
 & =2R_{f}(z)L_{N}^{\left(0\right)}R_{f}(z)+O\left(\frac{1}{N^{2}}\right).
\label{LL1}
\end{align}
Direct calculation gives that for any function $g_{k,j}(z)=z_{k}z_{j},k<j$,
on $H_{N}$
\begin{equation}
L_{N}^{\left(1\right)}g_{k,j}(z)=-\frac{\mu
 }{N}(z_{k}z_{j}(1-z_{k})+z_{k}z_{j}(1-z_{j}))+\frac{\mu
 }{N^{2}}z_{j}(1-z_{k}).
\end{equation}
 Using this formula we get that for any $z\in H_{N}$
\begin{align}
L_{N}^{\left(1\right)}R_{f}^{2}(z) & =-2R_{f}(z)\frac{\mu }{N^{2}}\sum \limits
 _{k}f(k/N)z_{k}(1-z_{k})\nonumber \\
 & +\frac{2\mu }{N^{4}}\sum \limits _{k<j}f(k/N)f(j/N)z_{j}(1-z_{k})+
\frac{\mu }{N^{4}}\sum \limits _{k}f^{2}(k/N)z_{k}(1-z_{k})\nonumber \\
 & =2R_{f}(z)L_{N}^{\left(1\right)}R_{f}(z)+O\left(\frac{1}{N^{2}}\right).
\label{LL2}
\end{align}

Summing the formulas (\ref{LL1}) and (\ref{LL2}) we obtain
\begin{equation}
\label{sum1}
L_{N}R_{f}^{2}(z)=2R_{f}(z)L_{N}R_{f}(z)+O\left(\frac{1}{N^{2}}\right).
\end{equation}
Lemma \ref{lb1} is proved.

\subsection{Tightness}

We will make all considerations
for the part ${\bf (i)}$ ($\alpha\neq\beta$ and $\mu=\mu/N$) of the
theorem.
All reasonings and conclusions are valid  for the part ${\bf
(ii)}$ ($\alpha=\beta$ and $\mu=\mu/N^2$) of the theorem
with evident changes.
We will denote by $P_{N,\lambda}^{\left(T\right)}$ the probability
distribution on the path space $D([0,T], H_{N})$ corresponding to
the process $\xi_{N}(t)$ and by $E_{N,\lambda}^{\left(T\right)}$ the
expectation with respect to this measure.

Theorem 4.1 in \cite{Mitoma} (see also Section~\ref{mit}) yields that tightness of the sequence
of $\{Q_{N, \lambda}^{\left(T\right)}\}_{N=2}^{\infty}$ will be
proved if we prove the same for a sequence of distributions of
one-dimensional projection $\{(\zeta_{N}(tN),f),t\in
[0,T]\}_{N=2}^{\infty}$  for every $f\in S({\textbf{R}})$. So fix
$f\in S({\textbf {R}})$ and consider the sequence of distributions
of the processes $(\zeta _{N}(tN), f),\, t\in [0,T]$. Note that the
probability distribution of a process $(\zeta _{N}(tN), f)$ is a
probability measure on $D([0,T],{\textbf {R}})$ the Skorokhod space
of real-valued functions.

By definition of the process $\zeta _{N}(tN)$ we have that
\[
(\zeta _{N}(tN),f)=\sum_{k}\xi_{N,k}(tN)
\int\limits_{k/N}^{(k+1)/N}f(x)dx.
\]
It is easy to see  that
$$
(\zeta_{N}(tN),f)=R_{f}(\xi_{N}(tN)) + \phi_{N}(t),
$$
where the random process $\phi_{N}(t)$ is bounded
$|\phi_{N}(t)|<C(f)/N$ for any  $t\geq 0$ and sufficiently large
$N$. Therefore it suffices to prove  tightness of the sequence of
distributions of random processes $\{R_{f}(\xi_{N}(tN)),\,t\in
[0,T]\}_{N=2}^{\infty}$.

Introduce two random processes
\begin{equation}
\label{mart1}
W_{f,N}(t)=R_{f}(\xi _{N}(tN))-R_{f}(\xi _{N}(0))-
N\int \limits_{0}^{t}L_{N}R_{f}(\xi _{N}(sN))ds,
\end{equation}
and
\[
V_{f,N}(t)=(W_{f,N}(t))^{2}-\int\limits_{0}^{t}Z_{f,N}(s)ds,
\]
where
\begin{equation}
\label{zz}
Z_{f,N}(s)=N\left(L_{N}R_{f}^{2}(\xi _{N}(sN))-2R_{f}(\xi
_{N}(sN))L_{N}R_{f}(\xi _{N}(sN))\right).
\end{equation}
It is well known (Theorem 2.6.3 in \cite{PreDeMas} or Lemma A1.5.1 in
\cite{KipLan}) that the processes $W_{f,N}(t)$ and $V_{f,N}(t)$ are
martingales.

The bound (\ref{b1}) yields that
\[
\left|N\int \limits _{\tau }^{\tau +\theta }L_{N}R_{f}(\xi
 _{N}(sN))ds\right|\leq C\theta \quad (a.s.)
\]
for any time moment $\tau $. Thus the sequence of probability
distributions of random processes
$\{N\int_{0}^{t}L_{N}R_{f}(\xi_{N}(sN))ds,\,t\in [0,T] \}_{N=2}^{\infty}$
 is tight  by Theorems~\ref{bil} and~\ref{ald} from Appendix.

The bound (\ref{b2}) yields that
\begin{equation}
\label{kvad}
E_{N,\lambda}^{\left(T\right)}(W_{f,N}(\tau +\theta )
-W_{f,N}(\tau ))^{2}=E_{N, \lambda}^{\left(T\right)}
\left(\int\limits_{\tau}^{\tau +
\theta }Z_{f,N}(s)ds\right)\leq \frac{C\theta }{N}.
\end{equation}
for any stopping time $\tau \geq 0$ since $V_{f,N}(s)$ is martingale.
Using this estimate and Chebyshev inequality we obtain that
the sequence of probability distributions
of  martingales  $\{W_{f,N}(t),\,t\in [0,T]\}_{N=2}^{\infty }$ is also tight
by  Theorem~\ref{ald}.
Thus the sequence of probability distributions
of the processes $\{R_{f}(\xi_{N}(tN)),\,t\in
[0,T]\}_{N=2}^{\infty}$ is tight by the equation (\ref{mart1})
and the assumption (\ref{initial}) and, hence, the sequence of probability
measures  $\{Q_{N, \lambda}^{\left(T\right)}\}_{N=2}^{\infty}$ is tight
 by Theorem 4.1 in \cite{Mitoma}.

\subsection{Characterization of  a limit point}

We are going to show now that  there is  a unique limit point
of the sequence  $\{Q_{N, \lambda}^{\left(T\right)}\}_{N=2}^{\infty}$
and this limit point is supported by trajectories
which are weak solutions of the partial differential equation
(\ref{eq:ur-asym}) in the sense of Definition \ref{weak}.

Let $f(s,x)\in C_{0,T}^{\infty}$ and denote
$$
R_{f}(t,\xi_{N}(tN))=\frac{1}{N}\sum_{k}\xi_{N,k}(tN)f(t,k/N),
$$
Define as before two random processes
$$
W_{f,N}'(t)=R_{f}(t,\xi _{N}(tN))-R_{f}(0,\xi _{N}(0))-
\int \limits_{0}^{t}\left(\partial/\partial s+NL_{N}\right)
R_{f}(s,\xi_{N}(sN))ds,
$$
and
$$
V_{f,N}'(t)=(W_{f,N}'(t))^{2}-\int \limits _{0}^{t}Z_{f,N}'(s)ds,
$$
where
$$
Z_{f,N}'(s)=N\left(L_{N}R_{f}^{2}(s,\xi _{N}(sN))-2R_{f}(s,\xi
_{N}(sN))L_{N}R_{f}(s,\xi _{N}(sN))\right).
$$
By Lemma A1.5.1 in \cite{KipLan} the processes $W_{f,N}'(t)$ and $V_{f,N}'(t)$
are  martingales.
It is easy to see that
\begin{multline}
\label{rasn1}
W_{f,N}'(t) =(\zeta (tN),f)- \int\limits_{0}^{t}(\zeta_{N}(sN),f_{s}+\lambda f_{x}+\mu f)ds - (\zeta (0),f)\\+
\mu\int\limits_{0}^{t}R_{f}(s,\xi^{2}(sN))ds+O\left(\frac{1}{N}\right).
\end{multline}
We are going to approximate the nonlinear term in (\ref{rasn1}) by some
quantities making sense in the space of generalised functions
since we treat the processes distributions
 as  probability measures on a space
$D([0,T], S'(\textbf{R}))$.
Let $\varkappa \in C_{0}^{\infty}(\R)$ be a non-negative function such that
$\int _{R}\varkappa (y)\, dy=1$.
Denote $\ke (y)= \varkappa \left(y/\varepsilon \right)/\varepsilon$,
for $0<\varepsilon\leq 1$ and let
$
(\ke \ast\varphi (s))(x)=\int_{{\textbf{R}}} \ke (x-y)\varphi (y,s)dy
$
be a convolution of a generalised function $\varphi (s,\cdot)$ with the test function $\ke (y)$.
\begin{lm}
\label{cut0}
The following uniform estimate holds
$$
\left|R_{f}(s,\xi ^{2}(sN))-((\varkappa _{\varepsilon }*\zeta
  _{N}(sN))^{2},f)\right|\leq F_{1}(\varepsilon )+
  F_{2}(\varepsilon N)
$$
where  the functions $F_1$ and $F_2$ do not depend on  $\xi$ and $s$
and
\[
  \lim _{\varepsilon \downarrow 0}F_{1}(\varepsilon )=
  \lim _{r\rightarrow+\infty}
  F_{2}(r)=0.
  \]
\end{lm}
\emph{Proof.}
For definiteness we assume
 that $\varkappa(x)=0$ for $x\in (-\infty,-1-\delta')
\cup (1+\delta',\infty)$ for some positive $\delta'$.
It is easy to see that if $x\in [k/N,(k+1)/N)$ for
some $k$, then
\begin{align*}
(\ke\ast \zeta_{N}(sN))(x)&=
\sum_{j}\xi_{N,j}(sN)\int\limits_{j/N}^{(j+1)/N}\ke(x-y)dy\\
&=\frac{1}{N}\sum\limits_{j}\ke ((k-j)/N)\xi_{N,j}(sN)+
g_1(N,\varepsilon,x,\xi_{N}(sN))
\end{align*}
where the function $g_1(N,\varepsilon,x,\xi_{N}(sN))$ can be bounded as follows
\begin{eqnarray*}
|g_1(N,\varepsilon,x,s,\xi_{N}(sN))| & \leq  & \sum _{j}\int \limits _{j/N}^{(j+1)/N}\left|\ke (x-y)-\ke \left(\frac{k-j}{N}\right)\right|\, dy\\
 & \leq  & \frac{1}{N}\sum _{m}\, 2\max _{w\in \left[(m-1)/N,m/N\right]}\frac{1}{\varepsilon ^{2}}\left|\varkappa '\left(\frac{w}{\varepsilon }\right)\right|\cdot \frac{1}{N}\\
 & = & \frac{2}{\left(N\varepsilon \right)^{2}}\sum _{m}
 \max _{v\in\left[(m-1)/(N\varepsilon ),m/(N\varepsilon )\right]}
 \left|\varkappa'\left(v\right)\right|
 \, =\,
 \frac{2\,U_{N\varepsilon }^{+}\left(|\varkappa '|\right)}
 {N\varepsilon }
\end{eqnarray*}
(see the (\ref{e-uds}) for the notation~$U^+$).
Note that if $\varepsilon$ is fixed then
$U_{N\varepsilon }^{+}\left(|\varkappa '|\right)=O(1)$
as $N\rightarrow\infty$.

This representation implies that
$$((\ke\ast \zeta_{N}(sN))^2,f)=
\frac{1}{N}\sum\limits_{k}f(k/N)\left(
\frac{1}{N}\sum\limits_{j}\ke((k-j)/N)\xi_{N,j}(sN)\right)^2
+O\left(\frac{1}{\eps N}\right)$$
Therefore
$$
|R_{f}(s,\xi^{2}(sN))-
((\ke\ast \zeta_{N}(sN))^2,f)|\leq J_{f,s}(\delta',\eps,N)
+
K(\varepsilon ,N)+O\left(\frac{1}{\eps N}\right),
$$
where
$$
J_{f,s}(\delta',\eps,N)=
\frac{2}{N}\sum\limits_{k}
|f(s,k/N)|\frac{1}{N}\sum\limits_{j:|j-k|<(1+\delta')\eps N}
\ke ((k-j)/N)|\xi_{N,k}(sN)-\xi_{N,j}(sN)|
$$
and
\[
K(\varepsilon ,N)=Const\cdot \left|
\frac{1}{N}\sum _{m}\ke \left(\frac{m}{N}\right)-1
\right|.
\]
Evidently,
$\ds K(\varepsilon ,N)=Const\cdot \left|
\frac{1}{N\varepsilon}\sum _{m}\varkappa \left(\frac{m}{N\varepsilon}
\right)-1
\right|$  and, hence,   $K(\varepsilon ,N)$ tends to $0$ in the limit
"$\varepsilon$ is fixed, $N\to\infty$".
So we can include   $K(\varepsilon ,N)$
into    $F_{2}(\varepsilon N)$.

Consider now the term $J_{f,s}(\delta',\eps,N)$.
Using monotonicity of trajectories we obtain
that for any $j$ such that $|j-k|<(1+\delta')\eps N$
$$
\left|\xi_{N,k}(sN)-\xi_{N,j}(sN)\right|\leq
\xi_{N,k-[(1+\delta')\eps N]}(sN)-\xi_{N,k+[(1+\delta')\eps N]}(sN).
$$
Thus we have that
$$
J_{f,s}(\delta',\eps,N)\leq
\frac{2}{N}\sum\limits_{k}
|f(s,k/N)|(\xi_{N,k-[(1+\delta')\eps N]}(sN)-\xi_{N,k+[(1+\delta')\eps N]}(sN)).
$$
Integrating by parts we get the following bound
\begin{align*}
J_{f,s}(\delta',\eps,N)&\leq \frac{2}{N}\sum\limits_{k}
\left(|f(s,(k+[(1+\delta')\eps N])/N)|-
|f(s,(k-[(1+\delta')\eps N])/N)|\right)\xi_{N,k}(sN)\\
& \leq\frac{2}{N}\sum\limits_{k}
|f(s,(k+[(1+\delta')\eps N])/N)-f(s,(k-[(1+\delta')\eps N])/N)| \\
& \leq Const\cdot  MD  \cdot (1+\delta')\varepsilon,
\end{align*}
where $M=\max _{x,s}\left|f_{x}(s,x)\right|$, $D$ is a  diameter of
$\mbox{supp}\,f_{x}(s,x)$.

Note that the last inequality is uniform in trajectories.
Lemma \ref{cut0} is proved.

\begin{lm}
\label{cut}
 For every $\delta>0$
$$
\limsup\limits_{\eps\ri 0}\limsup\limits_{N\ri \infty}
P_{\lambda,N}^{\left(T\right)}
\left(\left|\int\limits_{0}^{T}\Bigl(R_{f}(s,\xi^{2}(sN))-
((\ke\ast \zeta_{N}(sN))^2,f)\Bigr)\,ds\right|>\delta\right)=0.
$$
\end{lm}

The proof of this lemma is omitted because it is a direct consequence
of Lemma~\ref{cut0}.

It is easy to see that for any
$f\in C_{0,T}^{\infty}$
a map $F_{f,T,\eps}(\varphi ):D([0,T],S')\ri {\textbf{R}}_{+}$ defined by
$$
F_{f,T,\eps}(\varphi )=\left|
\int\limits_{0}^{T}\left(\varphi(s),f_{s}+\lambda f_{x}+\mu f)
-\mu ((\ke\ast\varphi (s))^2,f)\right)ds
+(\varphi (0),f)\right|
$$
is continuous, therefore for any $\delta>0$ the set  $\left\{ \varphi \in
D([0,T],S'):F_{f,T,\eps}(\varphi )>\delta \right\} $ is open and hence
\[
\limsup\limits_{\eps\ri 0}Q_{\lambda}^{\left(T\right)}\left(\varphi: F_{f,T,\eps}(\varphi )>\delta\right)\leq
\limsup\limits_{\eps\ri 0}\liminf_{N\ri \infty }Q_{N, \lambda}^{\left(T\right)}\left(\varphi: F_{f,T,\eps}(\varphi)>\delta \right),
\]
where $Q_{\lambda}^{\left(T\right)}$ is a limit point of the sequence
$\{Q_{N, \lambda}^{\left(T\right)}\}_{N=2}^{\infty }$.
Obviously that
\begin{align}
\label{qpp}
Q_{N, \lambda}^{\left(T\right)}\left(\varphi: F_{f,T,\eps}(\varphi)>\delta \right)& \leq P_{N, \lambda}^{\left(T\right)}
\left(\sup_{t\leq T}|W^{'}_{f,N}(t)| >\delta/2 \right)\\
&+P_{N, \lambda}^{\left(T\right)}\left(
\left|\int\limits_{0}^{T}(R_{f}(s,\xi^{2}(sN))-
((\ke\ast \zeta_{N}(sN))^2,f)ds\right|>\delta/2\right).\nonumber
\end{align}
It is easy to see that the bound (\ref{b2}) obtained in Lemma \ref{lb1}
for the process $Z_{f,N}(s)$ is  also valid for the process
$Z_{f,N}'(s)$, therefore for any $t$
$$
E_{N, \lambda}^{\left(T\right)}(W_{f,N}'(t))^{2}=
E_{N, \lambda}^{\left(T\right)}\left(\int\limits_{0}^{t }Z_{f,N}'(s)ds\right)\leq \frac{Ct}{N},
$$
since $V_{f,N}'(t)$ is a martingale.
Kolmogorov inequality implies that for any $\delta >0$
\begin{equation}
P_{N, \lambda}^{\left(T\right)}\left(\sup _{t\leq T}|W_{f,N}'(t)|\geq \delta \right)\leq
\delta^{-2}E_{N, \lambda}^{\left(T\right)}(W_{f,N}'(T))^{2}=\delta^{-2}
E_{N, \lambda}^{\left(T\right)}\left(\int \limits_{0}^{T}Z_{f,N}'(s)ds\right)\leq  \frac{CT}{N\delta ^{2}}.
\label{kolm}
\end{equation}
The second term in (\ref{qpp}) vanishes to zero by Lemma
\ref{cut} as $N\ri \infty$ and $\eps\ri 0$.
Therefore for any $f\in C_{0,T}^{\infty}$ and
$\delta >0$
\begin{equation}\label{e-Qd}
\limsup\limits_{\eps\ri 0}
Q_{\lambda}^{\left(T\right)}
\left(\varphi: F_{f,T,\eps}(\varphi )>\delta\right)=0.
\end{equation}
Let us prove  that we can replace the convolution
in $F_{f,T,\eps}$ by its limit which is well defined with respect to the
measure  $Q_{\lambda}^{\left(T\right)}$.
First of all we note that for  any $C>0$
$B_{C}({\textbf{R}})$ the set of measurable functions $h$
such that $\|h\|_{\infty}\leq C$  is a closed subset of $S'({\textbf{R}})$ in
both strong and weak topology. Indeed, consider a sequence of
functions $g_{n}\in B_{C},\, n\geq 1$ and assume this sequence
converges in $S'$ to some tempered distribution $G\in S'$. We are
going to show that this generalized function is determinated by some
measurable function bounded by the same constant $C$. It is easy to
see that for every $n\geq 1$
$$
\left|\int\limits_{\textbf{R}} g_{n}(x)f(x)dx\right|\leq
C\|f\|_{L^1},\,\, f \in S,
$$
where $\|\cdot\|_{L^1}$ is a norm in $L^1$ the space of all
integrable functions. So the limit linear functional $G$ on $S$ is
also continuous in $L^1$-norm
$$|G(f)|\leq C\|f\|_{L^1},\,\,f \in S.$$
The space $S$ is a linear subspace of $L^1$ therefore by Hahn-Banach
Theorem (Theorem III-5, \cite{RS}) the linear functional $G$ can be
extended to a continuous linear functional $\tilde{G}$ on $L^{1}$
with the same norm and such that $\tilde{G}|_{S}=G$. Using  the
theorem about the general   form of a continuous linear functional
on $L^1$ (\cite{RS}) we obtain that
$$
G(f)=\int\limits_{\textbf{R}}g(x)f(x)dx,\,\,  f\in
L^{1},
$$
where $g$ is a measurable bounded function. Obviously that
$\|g\|_{\infty}\leq C$.

Obviously that for any $N\geq 2$ and fixed $t \in [0,T]$ we have that
$Q_{N, \lambda}^{\left(T\right)}\left(\varphi (t,\cdot)\in B_1 \right)=1$,
where $\varphi (t,\cdot)=\varphi(t)$ is a coordinate  variable on
$D([0,T],S')$.
Therefore  if some subsequence
$\{Q_{N', \lambda}^{\left(T\right)}\}$  of
$\{Q_{N, \lambda}^{\left(T\right)}\}_{N=2}^{\infty}$ converges weakly to a
limit point $Q_{\lambda}^{\left(T\right)}$, then for any fixed $t\in [0,T]$
\begin{equation}
\label{closed}
Q_{\lambda}^{\left(T\right)}\left(\varphi(t)\in B_{1}\right)\geq
\limsup_{N'\ri \infty}Q_{N', \lambda}^{\left(T\right)}
\left(\varphi(t) \in B_{1} \right)=1,
\end{equation}
since $B_{1}$ is closed. 
Next lemma gives an important property of the convolution on the set
of bounded functions $L^{\infty}$.
\begin{lm}\label{closed2}
Fix $\varphi\in L^{\infty}$. Then

1) for any   $0\leq s\leq T$
$$((\ke\ast \varphi (s))^2,f)- \int_{\R}\varphi^2(s,x)f(x)\, dx \to 0
\qquad (\varepsilon\to 0)
$$

2)   for any   $0\leq t\leq T$
$$\int\limits_0^t\biggl\{
\,
((\ke\ast \varphi (s))^2,f)- \int_{\R}\varphi^2(s,x)f(x)\, dx
\,\biggr\}\, ds
\to 0
\qquad (\varepsilon\to 0).
$$

\end{lm}

\medskip

\emph{Proof of Lemma~\ref{closed2}.}
\begin{eqnarray*}
\left|((\ke\ast \varphi (s))^2,f) -
\int_{\R}\varphi^2(s,x)f(x)\, dx\right|&\leq &
 \int\limits_{\R}\,\biggl|\biggl(\bigl((\ke\ast \varphi
 (s)\bigr)(x)\biggr)^2 - \varphi^2(s,x)\biggr| \, |f(x)|
 \,dx\\
 &\leq& 2 \|\varphi\|_{\infty} \int\limits_{\R}\,\biggl|\bigl((\ke\ast \varphi
 (s)\bigr)(x) - \varphi(s,x)\biggr| \, |f(x)|
 \,dx
 \end{eqnarray*}
To get the last inequality we used the identity $a^2-b^2= (a+b)(a-
b)$ and the fact that  $\|\ke\ast\varphi\|_{\infty}\leq
 \|\varphi\|_{\infty}$. To finish the proof it suffices
 to apply a well-known result about convergence of
$(\ke\ast\varphi)(s,\cdot)$ to $\varphi(s)$ in $L^1_{loc}$.  This
proves the first statement of the lemma. To get the second statement
we use again the boundedness of $ \varphi$ and $\ke\ast\varphi$ and
apply the Lebesgue theorem.
Lemma \ref{closed2} is proved.

\smallskip

On the set  $\supp Q_{\lambda}^{\left(T\right)}$ we can define a
functional
$$
F^0_{f,T}(\varphi)=
\int\limits _{0}^{T}(\varphi (s),f_{s}+\lambda f_{x}+\mu f)
-\mu (\varphi^{2}(s),f))ds+(\varphi (0),f).
$$
The equation (\ref{closed}) and Lemma \ref{closed2} yield that for any
$\varphi\in \supp Q_{\lambda}^{\left(T\right)}\,\,$ 
$F_{f,T,\varepsilon}(\varphi)\to F^0_{f,T}(\varphi)$ as $\varepsilon\to
0$. This implies that for any $\delta_1>0$
$$
Q_{\lambda}^{\left(T\right)}
\left(\varphi: |F_{f,T,\eps}(\varphi )-F^0_{f,T}(\varphi)| >\delta_1\right)
\rightarrow 0\qquad (\varepsilon\to 0).
$$
Combining this with (\ref{e-Qd}) we get that  a limit point
of the sequence
$\{Q_{N,\lambda}^{\left(T\right)}\}_{N=2}^{\infty }$ is
 concentrated on the trajectories $\varphi (t,\cdot),\,t\in [0,T],$
taking values in the set of regular bounded functions
and satisfying the following integral equation
\[
\int\limits _{0}^{T}(\varphi (s),f_{s}+\lambda f_{x}+\mu f)
-\mu (\varphi^{2}(s),f))ds+(\varphi (0),f)=0
\]
for any $f\in C_{0,T}^{\infty}$.
It means that each such a trajectory is a weak solution
of the equation (\ref{eq:ur-asym}) in the sense
of Definition \ref{weak}.

\subsection{Uniqueness of a weak solution} 
\label{edinstv}

\paragraph{The first order equation.}

Using the method  of Theorem~1 in the celebrated paper  \cite{Oleinik} 
of Oleinik we will show here that for any  measurable bounded initial 
function $\psi(x)$ there might be at most one weak solution of the 
equation  (\ref{eq:ur-asym}) in the sense of Definition \ref{weak} and 
no entropy condition is required.

Let $u(t,x)$ and $v(t,x)$ be  two weak solutions of the
equation (\ref{eq:ur-asym}) in the region $[0,T]\times \textbf{R}$
with the same initial condition $\psi$  (not necessarily from $H$).
Definition \ref{weak} implies that
\begin{multline}
\label{eqt0}
\int\limits_{0}^{T}\int\limits_{\textbf{R}}((u(t,x)-v(t,x))(f_{t}(t,x)
+\lambda f_{x}(t,x) +\mu (1-u(t,x)-v(t,x))f(t,x))dxdt=0,
\end{multline}
for any
$f\in C_{0,T}^{\infty}$.
Consider the following sequence of equations
\begin{equation}
\label{eqt1}
f_{t}(t,x)+\lambda f_{x}(t,x)+g_{n}(t,x)f(t,x)=F(t,x),
\end{equation}
with any  infinitely  differentiable  function $F$ equal to zero
outside of a certain bounded region, lying in the half-plane
$t\geq\delta_1>0$, where $\delta_1$ is an arbitrary small number.
The functions $g_{n}(t,x)$ are uniformly bounded for all $x,t, n\geq
1$ and converges in $L^{1}_{loc}$ to the function $g(t,x)=\mu (1-
u(t,x)-v(t,x))$ as $n \ri \infty$. The solution  $f_{n}(t,x)\in
C_{0,T}^{\infty}$ of the equation (\ref{eqt1}) is given by the
following formula (formula (2.8) in \cite{Oleinik})
$$
f_{n}(t,x)=\int\limits_{T}^{t}F(s, x+\lambda (s-t))
\exp\left\{\int\limits_{t}^{s}g_{n}(\tau,x+
\lambda (\tau-t))d\tau\right\}ds.
$$
Equation (\ref{eqt0}) yields that
\begin{equation}
\label{eqt2}
\int\limits_{\textbf{R}_{+}}\int\limits_{\textbf{R}}(u(t,x)-v(t,x))F(t,x)dxdt
=\int\limits_{\textbf{R}_{+}}\int\limits_{\textbf{R}}
(g(t,x)-g_{n}(t,x))f(t,x))dxdt.
\end{equation}
The right side of (\ref{eqt2}) is arbitrary small
for sufficiently large $n$ and, since the left side of (\ref{eqt2})
does not depend on $n$, so it is equal to zero. Therefore
$u=v$, since $F$ is arbitrary.

An \emph{existence} of a weak solution of the equation (\ref{eq:ur-asym})
follows from the general theory for  quasilinear
equations of the first order (for example, Theorem 8 in \cite{Oleinik}).
In the particular case of the equation  (\ref{eq:ur-asym})
it is possible to obtain an explicit formula for a weak solution.
First of all we note that the Cauchy problem
\begin{equation}
u_{t}(t,x)=-\lambda u_{x}(t,x)+\mu (u^2(t,x)-u(t,x)),
\qquad
u(0,x)=\psi (x),
\label{eq:1-por}
\end{equation}
has a unique classical solution, if
$\psi\in C^1(\mathbf{R})$  and  there is an explicit
formula for this solution.
Indeed, using  substitution $\uc(t,x)=u(t,x-\lambda t)$ we transform the
equation (\ref{eq:1-por}) into the equation
$$\uc _{t}(t,x)=-\mu \uc (t,x)(1-\uc (t,x)),\qquad   \uc (0,x)=\psi (x).$$
Considering $x$ as a parameter we obtain an ordinary differential equation
which is solvable and the solution is given by:
\begin{equation}
\label{eq-uo-ex}
\uc(t,x)=\frac{\psi (x)e^{-\mu t}}{1-\psi (x)+\psi (x)e^{-\mu t}}.
\end{equation}
So, if $\psi\in C^1(\mathbf{R})$,
then a unique classical solution
of the equation (\ref{eq:1-por}) is given by the following formula
\begin{equation}
\label{solut1}
u(t,x)=\frac{\psi (x+\lambda t)e^{-\mu t}}{1-\psi (x+\lambda t)+
\psi (x+\lambda t)e^{-\mu t}}.
\end{equation}

If we approximate any measurable bounded function $g $ in
$L_{loc}^{1}$ by a sequence of smooth functions  $\{g_{n},  n\geq
1\}$, then  the sequence of corresponding  weak solutions
$\{u_n(t,x), n\geq 1\}$, where $u_{n}(t,x)$ is  defined by the
formula (\ref{solut1}) with $\psi=g_n$,  converges in $L_{loc}^{1}$
to  the weak solution of the equation with initial condition $g$ by
Theorem 11 in \cite{Oleinik} or  Theorem 1 in \cite{Kruzhkov}. It is
easy to show by direct calculation that the $L_{loc}^{1}$--limit of
the sequence $\{u_n(t,x), n\geq 1\}$ is given by the same formula
(\ref{solut1}) with $\psi=g$.

The  formula (\ref{solut1}) yields that  if $\psi\in H(\textbf{R})$,
then $u(t,\cdot)\in H(\textbf{R})$ as a function of $x$ for any
fixed $t\geq 0$. If  a function $u(t,x)$ is a weak solution of the
equation (\ref{eq:1-por}), then this function  is differentiable at
a point $(t,x)$ iff the initial condition $\psi(y)$ is
differentiable at the point $y=x+\lambda t$.

\paragraph{KPP-equation.}

The equation (\ref{eq:ur-KPP}) is a quasilinear
parabolic equations of the second order and is a
particular case of the famous KPP-equation in~\cite{KoPP}.
It is known that there exists a
unique weak solution $u(t,x)$ of the problem
$$
u_{t}(t,x)=\g u_{xx}(t,x)+\mu  (u^2(t,x)-u(t,x)),
\qquad u(0,x)=\psi(x),$$
for any bounded measurable initial function $\psi$, this
solution is in fact a unique classical solution
and if $\psi\in H(\mathbf{R})$,  then
$u(t,\cdot)\in H(\textbf{R})$ as a function of $x$ for
any fixed $t$.  We refer to the paper  \cite{Oleinik}
for more details.

\section{System with  fixed number of particles}
\label{stability}

In this section we deal with the situation ``$N$ is fixed, $t\rightarrow\infty$''.

\subsection{Proof of Theorem \ref{ergodic}}

 Let $\si (y,w)$ be the rate of transition from
the state $y=(y_{1},\ldots ,y_{N})\in \Gamma $ to the state
$w=(w_{1},\ldots,w_{N})\in \Gamma $ for the Markov $y(t)$ chain. Define
\(
\si (y)=\sum\limits _{w\not =y}\si (y,w).
\)
   From definition of the particle system it follows that
\[
\sigma (y)=(\alpha +\beta )N+\frac{\mu _{N}}{N}\sum _{(i,j)}I_{\left\{
 y_{i}>y_{j}\right\} }.
 \]
 Since $\ds \sum _{(i,j)}I_{\left\{ y_{i}>y_{j}\right\} }\leq N(N-1)/2$
we have uniformly in $y\in \Gamma $
\begin{equation}
\lmin \, \leq \, \sigma (y)\, \leq \, \sigma _{N}^{*}
\label{eq:lambda-below-above}
 \end{equation}
with $\lmin=(\alpha +\beta )N$ and
$\quad \sigma _{N}^{*}=(\alpha +\beta )N+\mu_{N}(N-1)/2$.
A discrete time Markov chain \emph{}$\left\{ Y(n),\, n=0,1,\ldots
 \right\} $
on the state space $S$ with transition probabilities
\begin{equation}
\label{eq:def-p}
\p (y,w)\equiv \P \left\{ Y(n+1)=w\, |\, Y(n)=y\right\} =\begin{cases}
 \ds \frac{\sigma (y,w)}{\sigma (y)}\, , & y\not =w\\
 0\, , & y=w,\end{cases}
 \end{equation}
 is an embedded Markov chain of the continuous time
 Markov chain~$\left(y(t),\, t\geq 0\right)$.

Theorem \ref{ergodic} is a consequence of the following statement.
\begin{lm}
\label{Deblin}
The Markov chain $\left\{ Y(n),\, n=0,1,\ldots\right\}$
is irreducible, aperiodic and  satisfies
to the Doeblin condition: there exist $\varepsilon >0$, $m_{0}\in \mathbf{N}$
and finite set $A\subset \Gamma$ such that
\begin{equation}
\label{eq:deb}
\P \left\{ Y(m_{0})\in A\, |\, Y(0)=Y_{0}\right\} \, \geq \varepsilon,
\end{equation}
for any $Y_{0}\in S$. Therefore this  Markov chain
is ergodic (\cite{FayMM}).
\end{lm}
\emph{Proof of Lemma~\ref{Deblin}.} We are going to show that
condition~(\ref{eq:deb}) holds with $A=\left\{ (0,\ldots ,0)\right\} $,
$m_{0}=N$,
\[
\varepsilon =\left(\frac{\min \left(\alpha ,\beta ,\mu _{N}/N\right)}
{\sigma_{N}^{*}}\right)^{N}>0.
 \]
The transition probabilities of the Markov chain $\left\{ Y(n),\,
n=0,1,\ldots\right\}$ are uniformly bounded from below in the
following sense: if a pair of states $(z,v)$ is such that $\si
(z,v)>0$ (or, equivalently, $\p (z,v)>0$) then (\ref{eq:def-p})
implies that
\begin{equation}
\label{eq:uni-bound}
\p (z,v)>\min\left(\alpha ,\beta,\mu _{N}/N\right)/\sigma _{N}^{*}.
\end{equation}
So to prove~(\ref{eq:deb}) we need only to show
 that for any $y$ there exists a sequence of states
\begin{equation}\label{eq-cepochka}
 v^{0}=y,\quad
v^{1},\quad v^{2},\quad \ldots ,\quad v^{N}=(0,\ldots ,0)
\end{equation}
which can be subsequently visited by the Markov chain
$\left\{ Y(n),\,n=0,1,\ldots\right\}$. The last means that
i.e. $\p (v^{n-1},v^{n})>0$  for every $n=1,\ldots,N$ and hence
$$
\P \left\{ Y(N)\in A\, |\, Y(0)=y\right\} \, \geq\,
\prod_{n=1}^{N}\p (v^{n-1},v^{n}) \,  \geq \, \left(\frac{\min \left(\alpha,
\beta, \mu _{N}/N\right)}{\sigma_{N}^{*}}\right)^{N}
$$
as a consequence of the uniform bound (\ref{eq:uni-bound}).

To prove existence of the sequence (\ref{eq-cepochka}) let us assume
first that $y=(y_{1},\ldots,y_{N})\not =0$. Choose and fix some $r$
such that $y_{r}=\ds \max _{i}y_{i}>0$. Denote by
\[
n_{0}=\#\left\{ j:\, y_{j}=0\right\}
\]
the number of left-most particles. Let the right-most particle $y_{r}$
move $n_{0}$ steps to right:
\[
v^{n}-v^{n-1}=e_{r}^{(\sn )},\quad \quad  n=1,\ldots,n_{0}.
\]
This can be done by using of jumps to
the nearest right state. So $Y(n_{0})=v^{n_{0}}$ has exactly $n_{0}$
particles at $0$ and $N-n_{0}$ particles out of $0$. Denote by
$i_{n_{0}+1}<i_{n_{0}+2}<\cdots <i_{N}$ indices of particles with
$v_{j_{a}}^{n_{0}}>0,\,a=n_{0}+1,\ldots,N$. Let now the Markov
chain $Y$ transfer each of these particles to~$0$:
\[
v^{a}-v^{a-1}=-v_{i_{a}}^{n_{0}}e_{i_{a}}^{(\sn )},\quad \quad
 a=n_{0}+1,\ldots,N.
\]
It is possible due to transitions provided by the interaction.

To complete the proof we need to consider the case $y=(y_{1},\ldots ,y_{N})=0$.
It is quite easy:
\[
v^{1}=e_{1}^{(\sn )},\quad v^{2}=2e_{1}^{(\sn )},\quad v^{3}=3e_{1}^{(\sn
 )},\quad \ldots ,\quad v^{N-1}=(N-1)e_{1}^{(\sn )},\quad v^{N}=0.
 \]
Proof of the lemma is over.

Denote by $\py =\left(\py (y),\, y=(y_{1},\ldots ,y_{N})\in \Gamma \right)$
a unique stationary distribution of the Markov chain
$\left\{ Y(n),\, n=0,1,\ldots\right\}$.

{\it The proof of  Theorem \ref{ergodic}\/} is now easy.
First of all let us show that the uniform bound
$\La (y)\geq \lmin $
implies existence of a stationary distribution for the Markov chain $y(t)$.
Indeed, it is easy to check that if  $\py $ is the stationary distribution
of the embedded Markov chain $Y$ and
$Q$ is the infinitesimal matrix for the chain $y(t)$, then a
vector with positive components
 $s=\left(s(w),w\in S\right)$ defined as
 \[
s(w)=\frac{\py (w)}{\La (w)},
\]
satisfies to the equation $sQ=0$.
So for existence of a stationary distribution of the chain  $y(t)$
it is sufficient to show  $\ds \sum _{w\in S}s(w)<+\infty $.
It is easy to check the last condition:
\[
\sum _{w\in S}s(w)=\sum _{w\in S}\frac{\py (w)}{\La (w)}\leq
\frac{1}{\lmin }\sum _{w\in S}\py (w)=\frac{1}{\lmin }.
\]
Therefore the
continuous-time Markov chain $\left(y(t),\, t\geq 0\right)$
has a stationary distribution~$\pi =\left(\pi (y),\, y\in \Gamma \right)$
of the following form
\[
\pi (y)=\frac{\ds \frac{\py (y)}{\sigma (y)}}{\ds \sum _{w\in \Gamma }
\frac{\py (w)}{\sigma (w)}}\, .
\]
Denote $p_{yw}(t)=\P \left\{ y(t)=w\, |\, y(0)=0\right\} $.
The next step is to prove that the continuous-time Markov chain
$y(t)$ is ergodic. To do this we show that the following Doeblin condition
holds: {\it
 for some $j_{0}\in S$ there exists $h>0$ and $0<\delta <1$
such that $p_{ij_{0}}(h)\geq \delta $ for all $i\in S$.
}
It is well-known (\cite{FayMM})
that  this condition implies ergodicity and moreover
\[
|p_{ij}(t)-\pi (y)|\, \leq \, \left(1-\delta \right)^{[t/h]}.
\]

Let  $\tau_{k},\,k\geq 0,$ be the time of stay of the Markov chain
$y(t)$ in $k-$th consecutive state.
Condition on the sequence of the chain states  $y_{k},\,k\geq 0,$ 
the joint distribution of
the random variables $\tau _{k},\, k\geq 0,$
coincides with the joint distribution
of independent random variables exponentially distributed with parameters
$\La (y_k),\, k=0,1,\ldots,n$, so
the transition probabilities of the chain $y(t)$ are
\[
p_{yw}(t)=\sum _{n}\sum _{(y\rightarrow w)}\P \left\{ (y\rightarrow w)\right\}
 \int\limits_{\Delta _{t}^{n}}\quad
e^{-\La (w)(t-t_{n})}
\prod\limits_{k=1}^{n}\La (y_{k-1})e^{-\La (y_{k-1})(t_{k}-t_{k-1})}
\, dt_{1}\ldots dt_{n},
\]
where $n$ corresponds to the number of jumps of the chain $y$ during the time
interval $[0,t]$, the inner sum
 is taken over all trajectories
 $(y\rightarrow w)=\left\{ y=y_{0},y_{1},\ldots ,y_{n}=w\right\} $
with $n$ jumps, integration is taken over
$\Delta _{n}^{t}=\left\{ 0=t_{0}\leq t_1\leq \cdots \leq t_{n}\leq t\right\} $,
and
$$
\P \left\{ (y\rightarrow w)\right\} =\p (y,y_{1})\p (y_{1},y_{2})\cdots \p (y_{n-1},w)
$$
is a probability of the corresponding path for the embedded chain.
The equation (\ref{eq:lambda-below-above}) implies that
the integrand in $p_{yw}(t)$ is uniformly bounded from below by the expression
\[
\left(\lmin \right)^{n}
\exp (-\lmax t_{1})\cdots \exp (-\lmax (t_{n}-t_{n-1}))\exp (-\lmax (t-t_{n})),
\]
and, hence,
\begin{align*}
\int\limits_{\Delta _{t}^{n}}
e^{-\La (w)(t-t_{n})}
\prod\limits_{k=1}^{n}\La (y_{k-1})e^{-\La (y_{k-1})(t_{k}-t_{k-1})}
\, dt_{1}\ldots dt_{n}& \geq
\frac{\left(\lmin t\right)^{n}}{n!}e^{-\lmax t}\\
 & = \P \left\{ \Pi _{t}=n\right\} e^{-\left(\lmax -\lmin \right)t},
\end{align*}
where $\Pi _{t}$ is a Poisson process with parameter $\lmin $.
It provides us with a lower bound
for the transition probabilities of the time-continuous chain:
\[
p_{yz_{0}}(t)\geq \left(\sum _{n}p^{n}(y,z_{0})\P \left\{ \Pi _{t}=n\right\} \right)e^{-\left(\lmax -\lmin \right)t}.
\]
It is easy now to get a lower bound for probabilities $p^{n}(y,z_{0})$.
Fix some $z_{0}\in S$ and denote $\mxi =\py (z_{0})$.
It follows from \underbar{\it ergodicity} of the chain $Y$, that for any fixed
$z_{0}$
\[
\p ^{m}(y,z_{0})\geq \py (z_{0})/2=\mxi/2>0,
\]
for all $m\geq m_{1}=m_{1}(y)$.
 For a \underbar{\it Doeblin}
Markov chain we have more strong conclusion, namely, the above number $m_{1}$
does not depend on $y$. Let us fix such $m_1$ and show that the
continuous-time Markov chain $y(t)$ satisfies to the Doeblin condition. Indeed,
\begin{eqnarray*}
p_{yz_{0}}(t) & \geq  & \left(\sum _{n<m_{1}}p^{n}(y,z_{0})
\P \left\{ \Pi _{t}=n\right\} +\sum _{n\geq m_{1}}p^{n}(y,z_{0})
\P \left\{ \Pi _{t}=n\right\} \right)e^{-\left(\lmax -\lmin \right)t}\\
 & \geq  & \sum _{n\geq m_{1}}p^{n}(y,z_{0})\P \left\{ \Pi _{t}=n\right\}
e^{-\left(\lmax -\lmin \right)t}\\
 & \geq  & \frac{\mxi }{2}\P \left\{ \Pi _{t}\geq m_{1}\right\}
e^{-\left(\lmax -\lmin \right)t}.
\end{eqnarray*}
Hence, the Doeblin condition holds: we choose any  $z_{0}$
as  $j_{0}$, take a corresponding $m_{1}$, fix any
$h>0$ and put
\[
\delta =\frac{\mxi }{2}\,
\P \left\{ \Pi _{h}\geq m_{1}\right\} e^{-\left(\lmax -\lmin \right)h}.
\]
Proof of the theorem is over.

\subsection{Evolution of the center of mass}

Consider the following function on the state space $\mathbf{Z}^{N}$:
$m(x_{1},\ldots ,x_{N})=\left(x_{1}+\cdots +x_{N}\right)/N.$
So if each\, particle\, has the mass~$1$\, and\, $x_{1}(t),\ldots ,x_{N}(t)$
are positions of particles, then $m(x_{1}(t),\ldots ,x_{N}(t))$ is
the \emph{center of mass} of the system.
We are interested in evolution of $\E \, m(x_{1}(t),\ldots ,x_{N}(t))$.
A direct calculation gives that
\begin{eqnarray}
\left(G_{N}m\right)(x_{1},\ldots ,x_{N}) & = & \sum
_{i=1}^{N}\left(\frac{\alpha
 }{N}-\frac{\beta }{N}\right)+\sum _{i=1}^{N}\sum _{j\not
 =i}\left(-\frac{x_{i}-x_{j}}{N}\right)I_{\left\{ x_{i}>x_{j}\right\}
}\frac{\mu
 _{N}}{N}\nonumber\\
 & = & (\alpha -\beta )-\frac{\mu _{N}}{N^{2}}\sum
 _{i<j}\left|x_{i}-x_{j}\right|,\label{eq:g-m}
\end{eqnarray}
where we have used the following equalities
\begin{eqnarray*}
m\left(x\pm e_{i}^{\left(\sn \right)}\right)-m\left(x\right)&=&
\pm \frac{1}{N}, \\
\quad m\left(x-(x_{i}-x_{j})e_{i}^{\left(\sn
 \right)}\right)-m\left(x\right)&=&-\frac{x_{i}-x_{j}}{N},\\
\left(x_{i}-x_{j}\right)I_{\left\{ x_{i}>x_{j}\right\}
 }+\left(x_{j}-x_{i}\right)I_{\left\{ x_{j}>x_{i}\right\} }&=&
 \left|x_{i}-x_{j}\right|.
\end{eqnarray*}

Note that the summand
 \[
-\frac{\mu _{N}}{N^{2}}\sum _{i<j}\left|x_{i}-x_{j}\right|
\]
added by the interaction to the {}``free dynamics'' drift $(\alpha
-\beta )$
depends only on the relative disposition of particles.
So the  center of mass of the system moves with speed which
tends to the value
\[
(\alpha -\beta )-\mu _{N}\, \frac{N-1}{2N}\, \, \E _{\pi }\,
 \left|x_{1}-x_{2}\right|
\]
as $t$ goes to infinity. Here $\mathbf{E}_{\pi }\, \left|x_{1}-
x_{2}\right|$ is the mean distance between two particles calculated
with respect to the stationary measure $\pi$ of the Markov chain
$Y$.

Using  this fact and Theorem  \ref{ergodic} we can describe the \emph{long time
behavior of the particle system} in the initial coordinates $x'$s as follows.
Theorem \ref{ergodic} means that  the system of stochastic interacting particles
possesses some relative stability. In coordinates~$y$ the system  approaches
exponentially fast its equilibrium state. In the meantime  the particle system
considered as a~"single  body"  moves with an asymptotically constant speed.
The speed differs from the mean drift of the free particle motion
and this difference is  due to the interaction between the particles.

\section{On travelling waves and long
time evolution of solution of PDE}
\label{travel}
We deal here with partial differential equation in variables $(t,x)\in\R_+\times\R$.
\begin{dfn} Function  $w=w(x)$ is called a travelling wave solution of some PDE
 if there exists $v\in \mathbf{R}$ such that
the function
\(u(t,x)=w(x-vt)\)
is a solution of this PDE.
The number $v$ is speed of the wave~$w$.
\end{dfn}

\noindent
We are interested only in the travelling waves
having the following properties:
\textbf{U1)} $w(x)\in[0,1]$;
\textbf{U2)} $w(x)$ and $dw(x)/dx$ have limits as  $x\rightarrow \pm
\infty $ and, besides, $w(-\infty )=1$ and  $w(+\infty )=0$.
We identify two travelling waves $w_{1}(x)$ and $w_{2}(x)$
if $w_{1}(x)=w_{2}(x-c)$ for some $c$.

For any probabilistic solution $0\leq u(t,x)\leq 1$ 
we define a function $r(t)$ such that
\(
u(t,r(t))\equiv \frac{1}{2}.
\)
Let a function $w(x)$ be a travelling wave solution.
Without loss of generality we can assume that $w(0)=1/2$.

\begin{dfn}
A solution  $u(t,x)$  converges \emph{in form}
to the travelling wave  $w(x)$ if
\[
\lim\limits _{t\rightarrow +\infty }u(t,x+r(t))=w(x),
\]
uniformly on any finite interval.
The solution $u(t,x)$ converges \emph{in speed} to the travelling wave  $w(x)$
if there exists $r'(t)=dr(t)/dt$
and
\( \lim\limits _{t\rightarrow +\infty }r'(t)=v,\)
where $v$ is a speed of the travelling wave $w(x)$.
\end{dfn}

\paragraph{First order equation.}

It is easy to check  for the equation~(\ref{eq:ur-asym}) for every $v<\lambda $
there exists a unique (up to shift) travelling wave solution
having properties  U1--U2 and this travelling wave solution
is given by the following formula
$\ds
w_v(x)=\left(1+\exp \left(\frac{\mu }{\lambda
 -v}x\right)\right)^{-1}.
$

\begin{prop}
\label{trav1}
If for some $C>0,\nu >0$, the initial profile $\psi (x)\in H(\mathbf{R})$ of the
equation (\ref{eq:ur-asym}) has the following asymptotic behavior
$
1-\psi (x)\sim C\exp (\nu x),
$
as $x\rightarrow -\infty$,   then there exists $x_0\in\R$ such
that for every $x\in \mathbf{R}$
$$
|u(t,x)-w_{v}(x-x_0-vt)|\rightarrow 0,
$$
as $t\rightarrow \infty$, where $\ds v=\lambda -\mu/\nu$.
\end{prop}
The proof of Proposition \ref{trav1} is a direct calculation
based on the exact formula (\ref{solut1}).

The formula (\ref{eq-uo-ex}) yields  that
$\uc(t_{1},x)\geq \uc(t_{2},x)$ for any $t_{1}<t_{2}$
and this observation  immediately  implies the following statement:
for every $x$ \
$\uc(t,x)\rightarrow I_{\left\{ y:\, \psi (y)=1\right\} }(x),$
as $t\rightarrow \infty$. As a direct application of this property we can obtain the following
\begin{prop}
\label{eq:for}
Assume that the initial profile has the form $\psi (x)=I_{\left\{ y<b\right\} }(x)$
 for some $b\in
\mathbf{R}$. Then
$$|u(t,x)-I_{\left\{ y<b\right\} }(x-\lambda t)|\rightarrow 0,\qquad t\rightarrow \infty.$$
\end{prop}
 So the function
$w(x)=I_{\left\{ y\leq 0\right\} }(x)$
is a unique (up to shift) non-increasing continuous
from the right travelling wave corresponding to the maximal possible speed
$v=\lambda $. It is easy to see that
this function is a limiting case of~$w_v(x)$ as $v\rightarrow\lambda-0$.

\paragraph{Second order equation.}
The existence of travelling waves for parabolic
partial differential equations was a subject of studying
in many papers followed to the paper~\cite{KoPP}.
A review of many results can be found in~\cite{Volpert} (see
also \cite{Smol}) and for completeness of the text we mention
some of them.
Reformulating the well-known results (\cite{Volpert}) we obtain
that travelling waves of the equation~(\ref{eq:ur-KPP}) can move only from the right to the
left. It means that the speed of any travelling wave is negative and, moreover,
 is bounded away from~0.
\begin{prop}
For equation~(\ref{eq:ur-KPP}) for every
$v\leq v_{*}=-\sqrt{4\g \mu}$ there exist
and unique (up to shift)
travelling wave solution with speed~$v$.
There are no other travelling wave solutions satisfying
the conditions U1 and U2.
\end{prop}

If a function $f=f(x)$ is such that $f(x)\leq 1$, $f(x)\rightarrow 1$
as  $x\rightarrow -\infty $ and there exists a limit
\(
\varkappa =\lim\limits _{x\rightarrow -\infty }{x}^{-1}{\log (1-f(x))}>0,
\)
then the number
$\varkappa $ is called \emph{Lyapunov exponent of the function $f$
(at minus infinity).}
It is well known (\cite{Volpert}) that for the equation~(\ref{eq:ur-KPP})
a travelling  $w(x)$ with speed~$v$ has the following Lyapunov exponent
at minus infinity
\[
\varkappa (v)=\left(-v-\sqrt{v^{2}-4\g \mu }\right)/\left(2\gamma\right).\]
Hence we get that for the travelling wave with minimal in absolute
value speed
$v_{*}=-\sqrt{4\g \mu }$ the Lyapunov exponent is
\(
\varkappa (v_{*})=\sqrt{{\mu }/{\g }}.
\)

\begin{prop}[\cite{Volpert}] Assume that  an initial function $\psi (x)$
has a  Lyapunov exponent $\varkappa $. Then

a) if $\varkappa \geq \sqrt{\mu/\g}$ then the solution $u(t,x)$
of the problem~(\ref{eq:ur-KPP}) converges in form and in speed
to the travelling wave moving with the minimal
speed $v_{*}=-\sqrt{4\g \mu }$;

b) if $\varkappa <\sqrt{\mu/\g}$ then the solution
 $u(t,x)$ of the problem~(\ref{eq:ur-KPP})
converges in form and in speed to the travelling wave
with speed
\(  \ds
v=-\left(\gamma\varkappa +\frac{\mu }{\varkappa }\right),\)
or, in other words, $\varkappa (v)=\varkappa $.

\end{prop}

We see from the above analysis that both the first order 
PDE~(\ref{eq:ur-asym}) and the second order PDE~(\ref{eq:ur-KPP}) 
exhibit  similar long-time behavior of their solutions. This seems very 
natural if we recall from Theorem~\ref{hydro} that the both equations 
arise as  hydrodynamical approximations of the same stochastic particle 
system.

\appendix

\section{Appendix}

\subsection{Strong topology on the Skorokhod space}
\label{top}
Remind that Schwartz space $S(\R )$ is a Frechet
space (see~\cite{RS}). In the dual space $S^{\prime }(\R )$ of tempered distributions
there are at least two ways to define topology (both not metrizable):

1) \emph{weak topology} on $S^{\prime }(\R )$, where all functionals
\(
\left( \, \cdot \, ,\phi \right) ,\quad \phi \in S(\R )
\)
are continuous.

2) \emph{strong topology} ($\strtop$) on $S^{\prime }(\R )$, which is generated by the set
of seminorms
\[
\left\{ \rho _{A}(M)=\sup _{\phi \in A}\left|\left( M,\phi
\right) \right|\, :\, \, A\subset S(\R )\, -\, {bounded}\right\} .
\]
Below we shall consider $S^{\prime }(\R )$ as equipped with the strong topology.
The problem of introducing of the Skorokhod topology on the space 
$D_{T}(S^{\prime}):=D([0,T],S^{\prime}(\R ))$ was studied in~\cite{Mitoma} and~\cite{Jakub}. 
We follow these papers.  For each seminorm $\rho _{A}$ on $S'(\textbf{R})$ we define the
following  pseudometric on $D([0,T], S'(\textbf{R}))$
\[
d_{A}(y,z)=\inf _{\lambda \in \Lambda }\left\{ \sup _{t}\rho_A(y_{t}-z_{\lambda (t)})
+\sup _{t\not =s}\left|\log \frac{\lambda (t)-\lambda (s)}{t-s}\right|\right\}
,\quad y,z\in
D([0,T], S'(\textbf{R})),\]
where $\inf$ is taken over the set
 $\Lambda =\left\{ \lambda =\lambda (t),\, t\in [0,T]\right\} $
of all strictly increasing continuous mappings of $[0,T]$ onto itself.
Introducing on $D([0,T], S'(\textbf{R}))$ the projective limit topology of
$\{d_A(\cdot,\cdot)\}$ we get a completely regular topological space.

\subsection{Mitoma theorem}
\label{mit}

Let $\mathcal{B}_{D_{T}(S^{\prime})}$ be the corresponding Borel $\sigma $-algebra. 
Let $\left\{ P_{n}\right\} $ be a sequence of probability measures
on $\left(D_{T}(S^{\prime}),\mathcal{B}_{D_{T}(S^{\prime})}\right)$. For each
$\phi \in S(\R )$ consider a map $\mathcal{I} _{\phi }:\, y\in
D_{T}(S^{\prime})\rightarrow (y,\phi )\in D_{T}(\R )$. The following
result belongs to I.~Mitoma~\cite{Mitoma}.

\begin{thr}

Suppose that for any $\phi \in S(\R )$ the sequence $\left\{ P_{n}\mathcal{I}
_{\phi }^{-1}\right\} $ is tight in  $D_{T}(\R )$. Then the 
sequence~$\left\{ P_{n}\right\} $ itself is tight in~$D_{T}(S^{\prime})$.
\end{thr}

\subsection{Probability measures on the Skorokhod space: tightness}

Let $\left\{ \left(\xi _{t}^{n},t\in [0,T]\right)\right\} _{n\in \mathbf{N}}$
\ be a sequence of real random processes which trajectories are
right-continuous and admit left-hand limits for every $0<t\leq T$ . We will
consider $\xi ^{n}$ as random elements with values in the Skorokhod space 
$D_{T}(\R ):=D\left([0,T],\R ^{1}\right)$ with the standard topology. Denote 
$P_{T}^{n}$ the distribution of $\xi ^{n}$, defined on the measurable space 
$\left(D_{T}(\R ),\mathcal{B}\left(D_{T}(\R )\right)\right)$. The following
result can be found in~\cite{billingsley}.

\begin{thr} \label{bil}
~\ The sequence of probability measures $\left\{ P_{T}^{n}\right\} _{n\in
\mathbf{N}}$ is tight iff the following two conditions hold:

1) for any $\varepsilon >0$ there is $\, C(\varepsilon )>0$ such that
\[
\sup _{n}P_{T}^{n}\left(\sup _{0\leq t\leq T}\left|\xi
_{t}^{n}\right|>C(\varepsilon )\right)\leq \varepsilon \, ;
\]

2) for any $\varepsilon >0$
\[
\lim _{\gamma \rightarrow 0}\limsup _{n}P_{T}^{n}\left(\xi _{\cdot }:\,
w^{\prime }(\xi ;\gamma )>\varepsilon \right)=0\, ,
\]
where for any function $f:\, [0,T]\rightarrow \R $ and any $\gamma >0$ we
define
\[
w^{\prime }(f;\gamma )=\inf _{\left\{ t_{i}\right\} _{i=1}^{r}}\, \max
_{i<r}\, \sup _{t_{i}\leq s<t<t_{i+1}}\left|f(t)-f(s)\right|\, ,
\]
moreover the $\inf $ is over all partitions of the interval $[0,T]$ such
that
\[
0=t_{0}<t_{1}<\cdots <t_{r}=T,\qquad t_{i}-t_{i-1}>\gamma ,\quad i=1,\ldots
,r.
\]
\end{thr}

The following theorem is known as the sufficient condition of Aldous~\cite[Proposition 1.6]{KipLan}.

\begin{thr} \label{ald}
Condition 2) of the previous theorem follows from the following condition
\[
\forall \varepsilon >0\qquad \lim _{\gamma \rightarrow 0}\, \limsup _{n}\,
\sup _{\tau \in \mathcal{R}_{T},\, \theta \leq \gamma }P_{T}^{n}\left(\left|\xi _{\tau
+\theta }-\xi _{\tau }\right|>\varepsilon \right)=0\, ,
\]
where $\mathcal{R}_{T}$ \ is the set of Markov moments (stopping times) not
exceeding $T$ .
\end{thr}

\end{document}